\documentclass[a4paper,twoside]{article}
\usepackage[utf8]{inputenc}
\usepackage[english]{babel}
\usepackage{hyperref}
\usepackage{amsmath,amsthm}
\usepackage{amssymb,amsfonts}
\usepackage{graphicx}
\usepackage{color}
\usepackage{xcolor}
\usepackage{framed}
\usepackage{tikz}
\usepackage[font=small]{caption}

\colorlet{shadecolor}{blue!0}
\definecolor{Red}{rgb}{0,0,0.3}

\usepackage{lastpage}

\usepackage{fancyhdr}
\renewcommand{\epsilon}{\varepsilon}

\newtheorem{theorem}{\underline{\textrm Theorem}}

\newenvironment{theo}
  {\begin{shaded}\begin{theorem}}
  {\end{theorem}\end{shaded}}

\makeatletter
\DeclareFontFamily{U}{tipa}{}
\DeclareFontShape{U}{tipa}{m}{n}{<->tipa10}{}
\newcommand{\arc@char}{{\usefont{U}{tipa}{m}{n}\symbol{62}}}%

\newcommand{\arc}[1]{\mathpalette\arc@arc{#1}}

\newcommand{\arc@arc}[2]{%
  \sbox0{$\m@th#1#2$}%
  \vbox{
    \hbox{\resizebox{\wd0}{\height}{\arc@char}}
    \nointerlineskip
    \box0
  }%
}
\makeatother

\begin{document}
\setlength\parindent{0pt}\thispagestyle{empty}
\begin{center}{\vspace{-2cm}\Large \textbf{ Surprising identities for the hypergeometric $\phantom{}_4 F_3$  function}}\end{center}

\begin{minipage}{0.49\textwidth}\begin{center}
\textbf{Jacopo D'Aurizio}\\
Dipartimento di Matematica, Università di Pisa\\
\href{jacopo.daurizio@gmail.com}{jacopo.daurizio@gmail.com} \texttt{+393479835357}\\
$\phantom{}$
\end{center}
\end{minipage}
\begin{minipage}{0.49\textwidth}\begin{center}
\textbf{Sabino Di Trani}\\
Dipartimento di Matematica, Università di Firenze\\
\href{sabino.ditrani@unifi.it}{sabino.ditrani@unifi.it} \texttt{+393276158611}\\
Gruppo INDAM GnSAGA
$\phantom{}$
\end{center}
\end{minipage}

\rule{\textwidth}{1pt}
\section{Introduction}
A well-known solution to the Basel problem (asking for a closed form for $\zeta(2)$) is essentially due to Euler and goes through the following identity:
$$ \zeta(2)=\sum_{n\geq 1}\frac{1}{n^2} = 3\sum_{n\geq 1}\frac{1}{n^2\binom{2n}{n}}. $$
It can be proved in a variety of ways, for instance through creative telescoping or by exploiting (shifted) Legendre polynomials. 
It is interesting to point out that the RHS of the previous identity is directly related with the Taylor series of the squared arcsine, 
from which $\zeta(2)=\frac{\pi^2}{6}$. The purpose of our work is to provide a fair amount of generalizations of this approach, in order
to exhibit some rather surprising closed forms for values of the hypergeometric function $\phantom{}_4 F_3$, defined as
$$ \phantom{}_4 F_3 \left(a,b,c,d;e,f,g;z\right) = \sum_{n\geq 0}\frac{(a)_n (b)_n (c)_n (d)_n }{n! (e)_n (f)_n (g)_n}\,z^n,\qquad (a)_n = \frac{\Gamma(a+n)}{\Gamma(a)}. $$
A key identity for the following manipulations is the main property of the Beta function:
$$ \int_{0}^{1} x^{a-1} (1-x)^{b-1}\,dx = B(a,b) = \frac{\Gamma(a)\,\Gamma(b)}{\Gamma(a+b)} $$
holds \cite{Askey} for every $a,b$ with positive real part. We will prove, in particular, that
$$ \phantom{}_4F_3\left(1,1,1,\frac{3}{2};\frac{5}{2},\frac{5}{2},\frac{5}{2};1\right) = {\frac{27}{2}\left(7\,\zeta(3)+(3-2G)\,\pi-12\right)} $$
with $G$ being Catalan's constant, and that $\phantom{}_4 F_3\left(1,1,1,\frac{3}{2};\frac{4}{3},\frac{5}{3},2;\frac{4z}{27}\right)$ is given 
by a combination of a squared arctangent and a squared logarithm, just like $\zeta(2)$ is deeply related with $\arcsin^2(z)$. 
In the former case we will employ a convolution argument, in the latter a functional identity for the dilogarithm function.
The results presented here follow from an attempt to use the convolution technique to compute a closed form for
$$ \sum_{n\geq 1}\frac{\cos n}{n^2\binom{3n}{n}}, $$
Instead of tackling this sum with the convolution approach, we computed it by exploiting an unexpected connection between $\text{Li}_2$ and $\phantom{}_4 F_3$.

\section{The first surprising identity about \texorpdfstring{$_4 F_3$}{4F3}}
The problem of finding a closed form for
$$ \phantom{}_4F_3\left(1,1,1,\frac{3}{2};\frac{5}{2},\frac{5}{2},\frac{5}{2};1\right) $$
appeared on \href{http://math.stackexchange.com/q/2123298/44121}{math.stackexchange.com}\cite{MSE} in February 2017. We will give a complete outline of our original solution,\\ based on the following Lemmas:
\begin{equation*} \int_{0}^{\pi/2}\sin(x)^{2n+3}\,dx = \frac{4^{n}(2n+2)}{(2n+3)(2n+1)\binom{2n}{n}}\tag{1}\end{equation*}
\begin{equation*} \arcsin^2(x)=\frac{1}{2}\sum_{n\geq 0}\frac{4^{n+1} x^{2n+2}}{(2n+2)(2n+1)\binom{2n}{n}}\tag{2}\end{equation*} 
If we integrate both sides of $(2)$ we get:
\begin{equation*} -2x+2\sqrt{1-x^2}\arcsin(x)+x\arcsin^2(x) = \frac{1}{2}\sum_{n\geq 0}\frac{4^{n+1} x^{2n+3}}{(2n+3)(2n+2)(2n+1)\binom{2n}{n}}\tag{3}\end{equation*}
We just have to gain an extra $\frac{1}{(2n+3)}$ factor. For such a purpose, we divide both sides of $(3)$ by $x$ and perform termwise integration again, leading to:
\begin{equation*} -4x+2\sqrt{1-x^2}\arcsin(x)+x\arcsin^2(x)+2\int_{0}^{\arcsin(x)}\frac{u\cos^2(u)}{\sin(u)}\,du\\= \frac{1}{2}\sum_{n\geq 0}\frac{4^{n+1} x^{2n+3}}{(2n+3)^2(2n+2)(2n+1)\binom{2n}{n}}\tag{4}\end{equation*}
Now we evaluate both sides of $(4)$ at $x=\sin\theta$ and exploit $(1)$ to perform $\int_{0}^{\pi/2}(\ldots)\, d\theta$. That leads to:
\begin{equation*} \sum_{n\geq 0}\frac{16^n}{(2n+3)^3(2n+1)^2\binom{2n}{n}^2}=(\pi-4)+\int_{0}^{\pi/2}\int_{0}^{\theta}\frac{u\cos^2(u)}{\sin(u)}\,du\,d\theta \tag{5}\end{equation*}
and we are extremely close to a conclusion, since the last integral boils down to $\int_{0}^{\pi/2}\int_{0}^{\theta}\frac{u}{\sin u}\,du\,d\theta$, that is well-known. Since
$$ \int\frac{u\,du}{\sin u} = C + u\left(\log(1-e^{iu})-\log(1+e^{iu})\right) +i \left(\text{Li}_2(e^{-iu})-\text{Li}_2(e^{iu})\right) $$
we have $\int_{0}^{\pi/2}\int_{0}^{\theta}\frac{u}{\sin u}\,du\,d\theta = -\pi G+\frac{7}{2}\zeta(3)$ and 
\begin{theo}
\begin{eqnarray*}\phantom{}_4F_3\left(1,1,1,\frac{3}{2};\frac{5}{2},\frac{5}{2},\frac{5}{2};1\right)&=&27\sum_{n\geq 0}\frac{16^n}{(2n+3)^3 (2n+1)^2 \binom{2n}{n}^2}\\&=&{\frac{27}{2}\left(7\,\zeta(3)+(3-2G)\,\pi-12\right)}\end{eqnarray*}
\end{theo}
where $G$ is Catalan's constant, $\sum_{n\geq 0}\frac{(-1)^n}{(2n+1)^2}$. Remarkably, in $2017$ no widespread mathematical software (Mathematica, Maple, Mathcad, Sage, $\ldots$) were able to derive such 
closed form, despite the proof only relies on a sort of convolution argument/ Parseval's identity. The same applies to the following identities we are going to prove. Before 
starting the next section, it is worth mentioning that a simplified version of the above technique also shows that
$$ \sum_{n\geq 1}\frac{16^n}{(2n+1)^2 (2n+3)^2 \binom{2n}{n}^2}=\pi-3,\qquad \sum_{n\geq 1}\frac{16^n}{n^2(2n+1)^2\binom{2n}{n}^2}=4(\pi-3) $$
i.e. connects $\pi$ with $\phantom{}_3 F_2\left(1,2,2;\frac{7}{2},\frac{7}{2};1\right)$  and $\phantom{}_3 F_2\left(1,1,1;\frac{5}{2},\frac{5}{2};1\right)$. 
These identities can also be proved by exploiting the contiguity relations \cite{Bailey} for $\phantom{}_{3} F_2$, as Computer Algebra Systems correctly recognize.

\section{The second surprising identity about \texorpdfstring{$_4 F_3$}{4F3}}
The series we are going to tackle in this section will provide something similar to Clausen's formula \cite{Vidunas}
$$\phantom{}_3 F_2\left(2c-2s-1,2s,c-\frac{1}{2};2c-1,c;x\right)=\phantom{}_2 F_1\left(c-s-\frac{1}{2},s;c;x\right)^2 $$
but for a $\phantom{}_4 F_3$ hypergeometric function. In particular, we will start dealing with
$$ \sum_{n\geq 1}\frac{z^n}{n^2 \binom{3n}{n}} = \frac{z}{3}\cdot\phantom{}_4 F_3\left(1,1,1,\frac{3}{2};\frac{4}{3},\frac{5}{3},2;\frac{4z}{27}\right). $$
The forthcoming manipulations can be summarized as follows:
\begin{itemize}
 \item We may express $\binom{3n}{n}^{-1}$ in terms of the Beta function and consider the usual integral representation of $B(n,2n)$, $\int_{0}^{1}x^{n-1}(1-x)^{2n-1}\,dx$;
 \item Assuming that $z$ is sufficiently close to zero, we are allowed to exchange the integral $\int_{0}^{1}(\ldots)\,dx$ and the series $\sum_{n\geq 1}$, converting 
 the original series in an integral of the form $\int_{0}^{1}\frac{\log p_z(x)+\log p_z(1-x)}{x}\,dx$, where $p_z(x)$ is the polynomial $1-zx^2+zx^3$;
 \item We may recall that $$\int_{0}^{1}-\log\left(1-\frac{x}{\alpha}\right)\frac{dx}{x}=\text{Li}_2\left(\frac{1}{\alpha}\right),$$
 where the dilogarithm function $\text{Li}_2(z)$ is the analytic continuation of $\sum_{n\geq 1}\frac{z^n}{n^2}$;
 \item We may recall that the dilogarithm function, according to D.Zagier \cite{Zagier}, has a good sense of humour,\\ i.e. fulfills a lot of interesting functional identities. Among them 
       we are interested in
       $$ \text{Li}_2(1-x)+\text{Li}_2(1-x^{-1}) = -\frac{1}{2}\log^2(x) $$
       that is straightforward to prove by differentiating both sides and considering the particular case $x=1$.\\ Through the change of variable $x=\frac{w-1}{w}$ such reflection formula 
       takes the form
       $$ \text{Li}_2\left(\frac{1}{w}\right) + \text{Li}_2\left(\frac{1}{1-w}\right) = -\frac{1}{2}\log^2\left(\frac{w-1}{w}\right) $$
       that is extremely well-suited for our purposes; 
 \item Indeed, if we assume that the roots of $p_z(x)$ are $\alpha_z,\beta_z,\gamma_z$, by the previous steps we get that the original series just depends on the following sum of squared logarithms:
       $$ \log^2\left(1-\frac{1}{\alpha_z}\right)+\log^2\left(1-\frac{1}{\beta_z}\right)+\log^2\left(1-\frac{1}{\gamma_z}\right). $$
\end{itemize}
Here we proceed in full detail:
\begin{eqnarray*}
 \sum_{n\geq 1}\frac{z^n}{n^2 \binom{3n}{n}} = \sum_{n\geq 1}\frac{2z^n \Gamma(n)\Gamma(2n)}{\Gamma(3n+1)}&=&\frac{2}{3}\sum_{n\geq 1}\frac{z^n}{n}B(n,2n)\\
 &=& \frac{2}{3}\sum_{n\geq 1}\int_{0}^{1}\frac{z^{n}x^{2n-1}(1-x)^{n-1}}{n}\,dx\\ 
 &=& \frac{2}{3}\int_{0}^{1}\left(\sum_{n\geq 1}\frac{\left(z x^2(1-x)\right)^n}{n}\right)\frac{dx}{x(1-x)} \\
 &=& \frac{2}{3}\int_{0}^{1}\frac{-\log(1-zx^2+zx^3)}{x(1-x)}\,dx \\
 \left[\frac{1}{x(1-x)}=\frac{1}{x}+\frac{1}{1-x}\right]\qquad &=& \frac{2}{3}\int_{0}^{1}\frac{-\log p_z(x)-\log p_z(1-x)}{x}\,dx \\
 \left[\text{Li}_2\text{ reflection formula}\right]\qquad&=& -\frac{1}{3}\sum_{\xi\in\left\{\alpha_z,\beta_z,\gamma_z\right\}}\log^2\left(1-\frac{1}{\xi}\right).
\end{eqnarray*}
If we consider the case $z=\frac{1}{2}$ we get that the roots of $p_z(x)$ are given by $-1$ and $1\pm i$. In particular:
\begin{theo}
 $$\frac{1}{6}\phantom{}_4 F_3\left(1,1,1,\frac{3}{2};\frac{4}{3},\frac{5}{3},2;\frac{2}{27}\right)=\sum_{n\geq 1}\frac{1}{n^2 2^n \binom{3n}{n}}={\frac{\pi^2}{24}-\frac{\log^2(2)}{2}}. $$
\end{theo}
Due to the dilogarithm reflection formula $\text{Li}_2(z)+\text{Li}_2(1-z)=\frac{\pi^2}{6}-\log(x)\log(1-x)$ we also get the identity:
\begin{equation*} \frac{\pi^2}{24}=\sum_{n\geq 1}\frac{1}{n^2 2^n}\left(1-\frac{1}{\binom{3n}{n}}\right).\tag{6}\end{equation*}
In a similar fashion, if we replace $x$ with $\sqrt{\frac{1-x}{2}}$ in both sides of $(2)$ and perform $\int_{0}^{1}\left(\ldots\right)\,dx$ we get:
\begin{equation*} \frac{\pi}{2}-1=\sum_{n\geq 1}\frac{2^n}{(n+1) n^2\binom{2n}{n}}.\tag{7}\end{equation*}
If we consider the case $z=\frac{1}{12}$ we get that the roots of $p_z(x)$ are given by $-2$ and $\frac{3\pm i\sqrt{15}}{2}$. In particular:
\begin{theo}
 $$ \frac{1}{36}\cdot \phantom{}_4 F_3\left(1,1,1,\frac{3}{2};\frac{4}{3},\frac{5}{3},2;\frac{1}{81}\right)=\sum_{n\geq 1}\frac{1}{n^2 12^n \binom{3n}{n}}={\frac{2}{3}\arctan^2\left(\frac{\sqrt{15}}{9}\right)-\frac{1}{2}\log^2\left(\frac{3}{2}\right)}. $$
\end{theo}
If we take some real number $m>\frac{1}{2}$ (this constraint ensures we stay in the domain of analytic continuation of the dilogarithm function, that is 
crucial for exploiting the dilogarithm reflection formula as we did) and enforce $-m$ to be a root of $p_z(x)$, we get $z=\frac{1}{m^2+m^3}$ and 
\begin{theo}[\textbf{Thai identity}]
 
\begin{eqnarray*}\frac{1}{3m^2+3m^3}\cdot\phantom{}_4 F_3\left(1,1,1,\frac{3}{2};\frac{4}{3},\frac{5}{3},2;\frac{4}{27(m^2+m^3)}\right)&=&\sum_{n\geq 1}\frac{1}{n^2 (m^2+m^3)^n \binom{3n}{n}} \\&=& {\frac{2}{3}\arctan^2\sqrt{\frac{3m-1}{(m+1)(2m-1)^2}}-\frac{1}{2}\log^2\left(1+\frac{1}{m}\right)}\end{eqnarray*}
 
\end{theo}
Through a suitable change of variable this identity proves that $\phantom{}_4 F_3\left(1,1,1,\frac{3}{2};\frac{4}{3},\frac{5}{3},2;z\right)$ 
has a closed form in terms of a squared logarithm and a squared arctangent. For instance, in the limit case $m=\frac{1}{2}$ we get:

\begin{equation*} \frac{8}{9}\cdot\phantom{}_4 F_3\left(1,1,1,\frac{3}{2};\frac{4}{3},\frac{5}{3},2;\frac{32}{81}\right)=\sum_{n\geq 1}\frac{8^n}{n^2 3^n \binom{3n}{n}} = \frac{\pi^2}{6}-\frac{\log^2(3)}{2}. \tag{8}\end{equation*}
Theorem $\mathbf{4}$ might have some unexpected consequences. The involved series are fast-convergent, hence they can be used to provide accurate numerical approximations 
of (squared) logarithms or arctangents. Additionally the arithmetic structure of the general term is pretty simple, hence Theorem $\mathbf{4}$ might be useful for 
estimating the irrationality measure of (squared) logarithms or arctangents. Moreover, such identity for $\phantom{}_4 F_3$ might be an important adjunct to the Wilf-Zeilberger 
recursion method for evaluating series involving binomial coefficients, since it provides a starting point for the recursion in a non-trivial case.

\section{An extension to \texorpdfstring{$_5 F_4$}{5F4}}
The technique introduced in the previous section is flexible enough to be able to deal with
$$ \sum_{n\geq 1}\frac{z^n}{n^2 \binom{4n}{n}} = \frac{z}{4}\cdot\phantom{}_5 F_4\left(1, 1, 1, \frac{4}{3}, \frac{5}{3}; \frac{5}{4}, \frac{3}{2}, \frac{7}{4}, 2; \frac{27z}{256}\right)$$
with just a minor fix. By defining $p_z(x)$ as $1-zx^3+zx^4$ and denoting its roots as $\alpha_z,\beta_z,\gamma_z,\delta_z$ we have:  
\begin{eqnarray*}
 \sum_{n\geq 1}\frac{z^n}{n^2 \binom{4n}{n}} = \sum_{n\geq 1}\frac{3z^n \Gamma(n)\Gamma(3n)}{\Gamma(4n+1)}&=&\frac{3}{4}\sum_{n\geq 1}\frac{z^n}{n}B(n,3n)\\
 &=& \frac{3}{4}\sum_{n\geq 1}\int_{0}^{1}\frac{z^{n}x^{3n-1}(1-x)^{n-1}}{n}\,dx\\ 
 &=& \frac{3}{4}\int_{0}^{1}\left(\sum_{n\geq 1}\frac{\left(z x^3(1-x)\right)^n}{n}\right)\frac{dx}{x(1-x)} \\
 &=& \frac{3}{4}\int_{0}^{1}\frac{-\log(1-zx^3+zx^4)}{x(1-x)}\,dx \\
 \left[\frac{1}{x(1-x)}=\frac{1}{x}+\frac{1}{1-x}\right]\qquad &=& \frac{3}{4}\int_{0}^{1}\frac{-\log p_z(x)-\log p_z(1-x)}{x}\,dx \\
 \left[\text{Li}_2\text{ reflection formula}\right]\qquad&=& -\frac{3}{2}\sum_{\xi\in\left\{\alpha_z,\beta_z,\gamma_z,\delta_z\right\}}\log^2\left(1-\frac{1}{\xi}\right).
\end{eqnarray*}
By imposing that $p_z(x)$ vanishes at $\alpha_z= -m$ we get $z=-\frac{1}{m^3+m^4}$. In this case, however, the final expression does not simplify as nicely as before, 
involving a fair amount of cube roots. The cubic function in Theorem $4$ and the quartic function arising in the $\phantom{}_5 F_4$ case are related to 
the so-called \textit{pull-back transformations}, already studied by Mitsuo Kato in \cite{Kato}.

\section{Back to \texorpdfstring{$_3 F_2$}{3F2}}
Surprisingly, many Computer Algebra Systems cannot seem to manage
$$ \sum_{n\geq 1}\frac{z^n}{n\binom{3n}{n}} = \frac{z}{3}\cdot\phantom{}_3 F_2\left(1,1,\frac{3}{2};\frac{4}{3},\frac{5}{3};\frac{4z}{27}\right) $$
either, but our method takes care of such hypergeometric function too, giving:
$$ \sum_{n\geq 1}\frac{z^n}{n\binom{3n}{n}} = \frac{2}{3}\int_{0}^{1}\frac{xz}{1-z(1-x)x^2}\,dx. $$
In the particular case $z=\frac{1}{2}$ we have:
\begin{theo}
$$\frac{1}{6}\cdot\phantom{}_3 F_2\left(1,1,\frac{3}{2};\frac{4}{3},\frac{5}{3};\frac{2}{27}\right)=\sum_{n\geq 1}\frac{1}{n 2^n \binom{3n}{n}} = {\frac{\pi}{10}-\frac{\log 2}{5}}$$ 
\end{theo}
and in general, for any $m>\frac{1}{2}$ the previous hypergeometric function evaluated at $z=\frac{1}{m^2+m^3}$ can be expressed in terms of logarithms and arctangents only.

\section{Acknowledgments}
We would like to thank Ruangkhaw Chaokha, a young student that has addressed the authors' interest towards the study 
of the interactions between the dilogarithm reflection formulas and closed forms for some values of $\phantom{}_4 F_3$. The name chosen for Theorem $4$ is a tribute to her. 
We are also very grateful to the reviewers for the suggested improvements and the given references about the pull-back transformations.
\end{document}